\documentclass[11pt,letter]{article}
\usepackage{color,graphicx}
\usepackage{citesort}
\usepackage{subfigure}
\usepackage{amsmath,amssymb,amsfonts}

\newtheorem{theorem}{Theorem}

\newtheorem{lemma}[theorem]{Lemma}

\begin{document}

\title{Photoacoustic Tomography in a Rectangular Reflecting Cavity}

\author{L. Kunyansky, B. Holman, B.~T. Cox}

\maketitle

\begin{abstract}
Almost all known image reconstruction algorithms for
photoacoustic and thermoacoustic tomography assume that the
acoustic waves leave the region of interest
after a finite time. This assumption is reasonable if the
reflections from the detectors and surrounding surfaces
can be neglected or filtered out (for example, by time-gating).
However, when the object is surrounded by
acoustically hard detector arrays, and/or by additional acoustic
mirrors, the acoustic waves will undergo multiple reflections.
(In the absence of absorption they would bounce around in such a
reverberant cavity forever).
This disallows the use of the existing free-space
reconstruction techniques.
This paper proposes a fast iterative reconstruction algorithm for
measurements made at the walls of a rectangular reverberant cavity.
We prove the convergence of the iterations under a certain sufficient
condition, and demonstrate the effectiveness and efficiency
of the algorithm in numerical simulations.

\end{abstract}


\section*{Introduction}

Photoacoustic tomography (PAT) and the closely related modality
thermoacoustic tomography (TAT) \cite{KrugerPAT,KrugerTAT,Oraev94,WangCRC,Beard2011}
are based on the photoacoustic effect, in which an acoustic wave is generated
by the absorption of an electromagnetic (EM) pulse.
The distinction between PAT and TAT is that the former uses visible or
near-infrared light, while the latter uses energy in the microwave region.
There are several endogenous chromophores which absorb in these wavelength ranges,
the most important of which are oxy- and deoxy-hemoglobin, and externally administered,
molecularly-targeted, chromophores can be used as contrast agents.
These emerging modalities are therefore attracting considerable interest
for molecular and functional imaging applications in pre-clinical and clinical imaging.

When EM pulses at these wavelengths are sent into biological tissue, the EM
energy will be absorbed and then thermalised. When this happens sufficiently rapidly, there will be
simultaneous, small, increases in the temperature and pressure
in the regions where the energy was absorbed. As tissue is an elastic medium, the
pressure increase propagates as an acoustic (ultrasonic) pulse, and can
be detected as a time series at the boundary of the tissue.
The image of the initial pressure distribution can subsequently be
reconstructed from these measurements by solving a certain inverse problem.
Since acoustic waves propagate through soft
tissues with little absorption or scattering, PAT and TAT yield high-resolution images
related to the EM properties of the tissue.
Such images cannot be obtained by purely optical or electrical (or electromagnetic)
techniques such as, e.g. diffuse optical tomography or electrical impedance tomography.

In order to facilitate the acoustic measurements the object is usually
submerged in water or a gel with acoustic properties close to those of the
tissues.
One method of measuring the acoustic pressure is to use a small single-element
detector to record the time-dependent acoustic pressure at a point.
In order to obtain enough data for the reconstruction, the detector is
scanned across an imaginary acquisition surface surrounding the region of interest;
the EM pulse must be re-emitted for each new location of the detector.
The wave propagation within such an acquisition scheme can be modelled by
the free-space wave equation, as the reflection of the waves from the single
detector does not affect the measurement at that detector.
The corresponding inverse problem has been studied extensively (see, for example reviews
\cite{Handbook,WangCRC,KuKuREW} and references therein), and various
reconstruction techniques have been developed to solve it. These techniques include
time reversal \cite{HKN,XuWang04,burg-exac-appro}, series expansion methods
\cite{Norton1,Norton2,Kun-fast,Kun-ser,Haltm-circ,AK}, and several explicit
inversion formulas \cite{Finch04,Finch07,MXW2,Kunyansky,Kun-cube,nguyen} (the
references above are by no means exhaustive).

In order to significantly speed up the data acquisition, the
measurements may be made by an array of detectors. (Both piezoelectric and
optically-addressed arrays have been used). Such an array typically forms a
solid surface, a side effect of which is that the sound waves are
reflected back into the region of interest. In the case of a single planar array the
free space reconstruction methods still apply, since the
reflected waves will propagate away from the detector and will not affect the
measurements.
However, it is not possible to obtain exact reconstructions from
data measured on a finite section of a plane.
To overcome this limitation, more advanced measuring systems will have the object
surrounded by the detector arrays (or perhaps intentionally installed acoustic
mirrors, e.g.~\cite{Cox2007}). One such system, based on optically-addressed Fabry-Perot detection arrays,
is currently being investigated.
The waves within the reverberant cavity formed by the detectors and passive reflectors will undergo
multiple reflections. Under the idealized assumption of negligible acoustic absorption
used by most classical methods, the oscillations within such a cavity will
continue forever. All the standard reconstruction techniques assume that the
signal either has a final duration (due to Huygens' principle in 3D) or dies
out sufficiently fast (in 2D or in the presence of nonuniform speed of sound).
As a result, these techniques are not applicable in the presence of such
multiple reflections.

In the present paper we develop an algorithm for reconstructing the initial
pressure distribution within a rectangular reverberant cavity from acoustic pressure
time series measurements made on at least three mutually adjacent walls. The algorithm
finds first a crude initial approximation that can be further refined
iteratively. The derivation of the algorithm is presented in Section
\ref{S:derivation}, together with the convergence analysis for the iterations
(Section \ref{S:convergence}). The algorithm was tested in numerical
simulations (Section \ref{S:numerics}) that showed both the high quality of
the reconstruction and very low noise sensitivity of the proposed technique.

High resolution images in 3D utilized in modern biomedical practice are
described by hundreds millions of unknowns. In order to be useful in practice,
a reconstruction algorithm must be fast. The proposed technique is
asymptotically fast; it requires $\mathcal{O}(N^{3}\log N)$ floating point
operations (flops) per iteration on an $N\times N\times N$ computational grid.
In our simulations the reconstruction time was comparable with that of known
fast algorithms for various free-space problem (e.g. \cite{Kun-ser,Kun-fast}).

\section{Formulation of the problem}

\subsection{Photoacoustic Tomography in Free Space}

\label{S:freespace}
In conventional photoacoustic tomography, time domain measurements of the
photoacoustically generated acoustic waves are made by an array of ultrasound
detectors positioned on a measurement surface. The inverse problem of PAT/TAT
consists in finding the initial pressure
distribution from the measurements. Several
idealizing assumptions are typically made to simplify this problem. First, the
speed of sound is frequently assumed to be known and constant throughout the
whole space, i.e. the presence of any physical boundaries is neglected. Second, it
is assumed that the measurements are done at all points lying on some
imaginary observation surface surrounding the object of interest.
Third, the detectors are considered small and non-reflecting (acoustically transparent).
Under these assumptions the
acoustic waves can be considered as propagating in free space ($\mathbb{R}
^{n}$, $n=2$ or $3$) and the acoustic pressure $p=p(t,x)$ is a solution
to the following (free-space) initial value problem (IVP):
\begin{align}
\left(  \frac{\partial^{2}}{\partial t^{2}}-c_{0}^{2}\Delta\right)  p  &
=0,\quad x\in\mathbb{R}^{n},\quad t\in\lbrack0,\infty),\label{E:IVP}\\
p|_{t=0}  &  =f(x),\quad\left.  \frac{\partial p}{\partial t}\right\vert
_{t=0}=0 \label{E:IVP1}
\end{align}
where $f(x)$ is the initial acoustic pressure distribution. Given time series
measurements $g$ of the acoustic pressure at each point of the observation
surface $S$ surrounding the region of interest $\Omega$ (within which $f(x)$
is supported)
\begin{equation}
g=p(t,x)\quad x\in S,\quad t\in\lbrack0,T], \notag
\end{equation}
it is possible to reconstruct $f(x)$ exactly. Indeed, since $S$ is not a
physical boundary, for a sufficiently large time $T$, the pressure
$p$ in $\Omega$ will vanish. (In 3D, due to the Huygens principle, the
pressure will vanish after $T=\mathrm{diam}(\Omega)/c_{0}$. In 2D,
the pressure does not completely vanish in finite time; however, it will decay
fast enough that it can be approximately set to 0 at a sufficiently large
value of $T\geq\mathrm{diam}(\Omega)/c_{0}$.)

Now one can solve in $\Omega$ the following initial/boundary value problem
(IBVP)
\begin{align}
\left(  \frac{\partial^{2}}{\partial t^{2}}-c_{0}^{2}\Delta\right)  p  &
=0,\quad t\in\lbrack0,T],\quad x\in\Omega,\label{E:InverseIBVP1}\\
\left.  p\right\vert _{t=T}  &  =0,\quad\left.  \frac{\partial p}{\partial
t}\right\vert _{t=T}=0,\quad\left.  p\right\vert _{x\in S}=g(t,x),
\label{E:InverseIBVP2}%
\end{align}
backwards in time (from $t=T$ to $t=0)$, thus obtaining
$f(x)=p(0,x)$. This technique is called \textit{time reversal}. It is
theoretically exact and works for an arbitrary closed surface $S$; many of the
existing reconstruction methods are based on this idea. For instance, one can
solve IBVP (\ref{E:InverseIBVP1},\ref{E:InverseIBVP2}) directly using
numerical techniques~\cite{Finch04,XuWang04,HKN}. Eigenfunction decomposition
techniques~\cite{AK,Kun-ser} and some of the known inversion
formulas~\cite{Kun-cube} yield reconstructions that are theoretically
equivalent to the ones obtained by time reversal. In each case the crucial
underlying assumption is that the values of the pressure are recorded until
the acoustic wave vanishes within $\Omega$.

\subsection{Photoacoustic Tomography in a Reflective Cavity}

The methods of the previous section can be used so long as the detectors (or
anything else such as the tank walls) do not reflect the acoustic waves.
When an array of detectors is used, instead of a single scanned detector,
the situation will change as the array itself will reflect the wave.
In the case of a single planar array the traditional
reconstruction techniques still apply, since the reflected waves will
propagate away from the detector and will not affect the measurements.
However, the images reconstructed from measurements made using a single
finite-sized planar array configuration will contain `partial view' or
`limited data' artifacts because the acoustic waves travelling parallel (or
almost parallel) to the surface are not measured. In order to improve the
reconstruction one could either reflect those waves back on the detector array
using passive reflectors~\cite{Cox2007}, or add new (perpendicular) detector
arrays. If the region of interest is surrounded by the arrays and/or
reflectors, a reverberant cavity is formed. (Note that the free surface of a
liquid also reflects waves). Wave propagation in a reverberant cavity is no
longer represented by the IVP~(\ref{E:IVP}), and a different mathematical model is needed.

The proper model should take into account that the domain $\Omega$ is
surrounded by a reflecting boundary $\partial\Omega$, and the corresponding
boundary conditions should be imposed on the wave equation (as opposed to the
free space propagation discussed in the previous section). The measurement
surface $S$ is now a subset of the boundary $\partial\Omega$.

When the boundary is treated as sound-hard the acoustic pressure $p(t,x)$ is
 a solution to the following initial boundary value problem:
\begin{align}
\left(  \frac{\partial^{2}}{\partial t^{2}}-c_{0}^{2}\Delta\right)  p  &
=0,\quad x\in\Omega\label{E:Forward1}\\
p(0,x)  &  =f(x),\quad\frac{\partial p(0,x)}{\partial t}=0,\label{E:Forward2}
\\
\frac{\partial}{\partial n}p(t,x)  &  =0,\quad x\in\partial\Omega.
\label{E:Forward3}
\end{align}
Now the measurements can be written as
\begin{equation}
g=p(t,x)\quad x\in S\subseteq\partial\Omega,\quad t\in\lbrack0,T]. \notag
\end{equation}

Other boundary conditions might be appropriate in some circumstances, such as
for instance if the cavity were built as a tank and one part of the boundary
was a water-air interface where $p\approx0$, i.e. Dirichlet conditions are
imposed on this part of the boundary.
Also, it is not necessary to measure
the pressure on the whole boundary $\partial\Omega$, and so $S$ is only a
subset of $\partial\Omega$. (Taking this to an extreme, Cox et al.
\cite{Cox2009} proposed using a single measurement point when the cavity is
ray-chaotic. This scheme, however, is unlikely to yield a stable reconstruction.)

A distinct property of this model is the preservation of the acoustic energy
trapped within the cavity. Since the model assumes that there is no absorption
of the acoustic energy by the medium, and the (Neumann or Dirichlet) boundary
conditions correspond to a complete reflection of waves, the oscillations
within $\Omega$ will (theoretically) continue forever. In practice, of course,
this is not the case and the waves will soon decrease to the extent that
further measurements will become impossible due to the low signal-to-noise
ratio. While we are not explicitly modeling the absorption, we have to assume
that the measurement time $T$ is bounded. It will be chosen to correspond,
roughly, to several bounces of the acoustic wave between the cavity walls.

The preservation of acoustic energy within the reverberating cavity makes it
impossible to solve the inverse problem by the time reversal techniques
mentioned in Section~\ref{S:freespace}. Indeed, in order to obtain the
accurate reconstruction of $f$ one has to accurately prescribe conditions
$p(T,x)$, and $\partial p/\partial t(T,x)$ to initialize the time
reversal. However, there is no way to measure these data within the object.
Moreover, this values are of the same order of magnitude as $f$, and if one
simply replaces them by a zero (as we can safely do in the free space case)
the induced error will be also of the same magnitude as $f$ we seek to
reconstruct. Therefore, almost none of the known reconstruction algorithms are
applicable here (with the exception
of~\cite{Cox2007,Cox2009,Ammari,WangYang07}). Below we propose inversion
techniques suitable for photoacoustic tomography within the reverberant cavity.

\section{Derivation of the reconstruction algorithm\label{S:derivation}}

In this section we develop a fast and accurate reconstruction algorithm for
photoacoustic tomography within a rectangular reverberant cavity. On one
hand, such an acquisition geometry is one of the simplest from the
mathematical standpoint, allowing one to design a simple and fast Fourier
based reconstruction technique.
On the other hand, photoacoustic scanners with this particular configuration are quite conceivable -
indeed a design based on optically-addressed Fabry-Perot ultrasound arrays~\cite{Zhang08}
is currently under development - and algorithms will be needed to process the data.

The reconstruction technique we propose can be used in both 2D and 3D settings,
with only minor changes. The 3D case is, obviously, more interesting from the
applied point of view, while the 2D case is somewhat easier to present. Thus,
in the next section we provide a derivation for the 2D algorithm, with the
caveat that the 3D case is quite similar. In Section~\ref{S:numerics} we test the 3D
version of the algorithm in numerical simulations.

\subsection{2D formulation and the series solution of the forward problem}

\begin{figure}[t]
\begin{center}
\includegraphics[width=2.4in,height=2.0in]{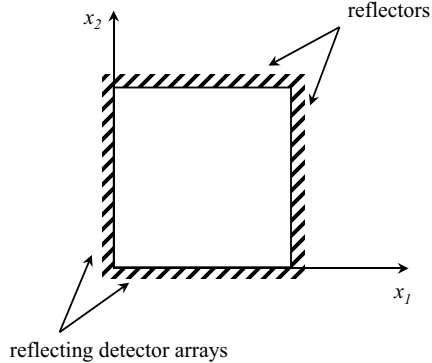}
\end{center}
\caption{
The 2D reverberant cavity with two passive reflecting walls, and two
reflecting detector arrays on which the data $g_{1}(t,x_2)$ and $g_{2}(t,x_1)$ are measured.}
\label{F:cavity}
\end{figure}

We assume that the acoustic pressure $p(t,x),$ $x=(x_{1},x_{2})$ solves the
initial/boundary value problem (\ref{E:Forward1})-(\ref{E:Forward3}) in 2D, in
a square $\Omega=[0,1]\times\lbrack0,1]$. Without loss of generality the speed
of sound $c_{0}$ will be assumed to equal 1. The measurements are represented
by the pair of time series $\mathbf{g}=(g_1(t,x_2),g_2(t,x_1))$
corresponding to the pressure values on the two adjacent
sides of the square $\Omega$ (Figure~\ref{F:cavity}):
\begin{equation}
g_{1}(t,x_{2})=\left.  p(t,x)\right\vert _{x_{1}=0},\qquad g_{2}%
(t,x_{1})=\left.  p(t,x)\right\vert _{x_{2}=0},\quad t\in\lbrack0,T].
\notag \end{equation}
Let us denote\ by $\mathcal{W}$ the linear operator transforming the initial
condition $f(x)$ into the boundary data $\mathbf{g}$:%
\begin{equation}
\mathbf{g}=\mathcal{W}(f).
\notag \end{equation}
Our ultimate goal is to reconstruct the initial pressure $f(x)=p(0,x)$
from $\mathbf{g},$ i.e. to invert $W$.

The simple shape of domain $\Omega$ allows us to utilize the eigenfunctions
\begin{equation}
\varphi_{k,l}(x)=\cos(\pi kx_{1})\cos(\pi lx_{2})
\notag \end{equation}
of the Neumann Laplacian in a square%
\begin{equation}
(\Delta+\omega_{k,l}^{2})\varphi_{k,l}=0,\quad x\in\Omega,\qquad\frac
{\partial\varphi_{k,l}}{\partial n}=0,\quad x\in\partial\Omega,
\notag \end{equation}
with the eigenfrequencies $\omega_{k,l}$:%
\begin{equation}
\omega_{k,l}=\pi\sqrt{k^{2}+l^{2}},\quad k,l=0,1,2,...
\notag \end{equation}

Using the eigenpairs $(\varphi_{k,l},\omega_{k,l}),$ the solution $p(t,x)$ of
the forward IBVP (\ref{E:Forward1})-(\ref{E:Forward3}) can be represented by
the following series:
\begin{equation}
p(t,x)=\sum_{k=0}^{\infty}\sum_{l=0}^{\infty}f_{k,l}\cos(\pi kx_{1})\cos(\pi
lx_{2})\cos(\omega_{k,l}t),\quad x\in\Omega,\quad t\in\lbrack0,T],
\label{E:series_sol}%
\end{equation}
where numbers $f_{k,l}$ are the coefficients of the two-dimensional Fourier
cosine series of the initial pressure $f(x)$:%
\begin{align}
f_{k,l}  &  =\frac{1}{||\varphi_{k,l}||_{2}^{2}}\int\limits_{\Omega
}f(x)\varphi_{k,l}(x)~dx,\quad k,l=0,1,2,...,\label{E:fourser1}\\
||\varphi_{k,l}||_{2}^{2}  &  =\int\limits_{\Omega}\varphi_{k,l}%
^{2}(x)~dx=\left\{
\begin{array}
[c]{cc}%
\frac{1}{4}, & k\neq0\text{ and }l\neq0,\\
1, & k=0\text{ and }l=0,\\
\frac{1}{2}, & \text{otherwise}%
\end{array}
\right.  . \label{E:fourser2}%
\end{align}
For future use we will denote the double-indexed sequence of the Fourier
coefficients by $F=\{f_{k,l}\}_{k,l=0}^{\infty},$ and will use the notation
$\mathcal{F}_{\operatorname{series}}$ to refer to the transformation defined
by equations (\ref{E:fourser1}) and (\ref{E:fourser2}), so that%
\begin{equation}
F=\mathcal{F}_{\operatorname{series}}f.
\notag \end{equation}

Obviously, in order to obtain $f(x)$ it is enough to reconstruct coefficients
$F;$ function $f$ is then computed by inverting $\mathcal{F}%
_{\operatorname{series}},$ i.e. by summing the Fourier cosine series:%
\begin{equation}
f(x)=\mathcal{F}_{\operatorname{series}}^{-1}F=\sum_{k=0}^{\infty}\sum
_{l=0}^{\infty}f_{k,l}\cos(\pi kx_{1})\cos(\pi lx_{2}).
\notag \end{equation}

\subsection{What can be found from the data measured on one side?}

Let us first analyze the connection between coefficients $f_{k,l}$ and the
data measured on one side of the square (say, $g_{2}$):
\begin{equation}
g_{2}(t,x_{1})\equiv\left.  p(t,x)\right\vert _{x_{2}=0}=\sum_{k=0}^{\infty
}\sum_{l=0}^{\infty}f_{k,l}\cos(\pi kx_{1})\cos(\omega_{k,l}t),\quad
t\in\lbrack0,T]. \label{E:oneside}%
\end{equation}
Due to the orthogonality of the cosine functions (in variable $x_{1})$ the
above equation splits:
\begin{align}
\sum_{l=0}^{\infty}f_{k,l}\cos(\omega_{k,l}t)  &  =g_{2,k}(t),\quad
t\in\lbrack0,T],\quad k=0,1,2,...,\label{E:NUFourier}\\
g_{2,0}(t)  &  \equiv\int\limits_{0}^{1}g_{2}(x_{1},t)dx_{1},\quad
g_{2,k}(t)\equiv2\int\limits_{0}^{1}g_{2}(x_{2},t)\cos(\pi kx_{1})dx_{1},\quad
k=1,2,3,.... \label{E:modes}%
\end{align}
It follows from equation\ (\ref{E:NUFourier}) that each of the $k$ functions
$g_{2,k}(t)$ is a combination of cosine functions in $t$ with known
frequencies; the values $f_{k,l}$ can be viewed as generalized Fourier
coefficients. More precisely, equation\ (\ref{E:NUFourier}) can be re-written
as the inverse Fourier transform $\mathcal{F}^{-1}$ of a sequence of Dirac
delta-functions:
\begin{equation}
g_{2,k}(t)=\sqrt{2\pi}\mathcal{F}^{-1}\left[  \sum_{l=0}^{\infty}f_{k,l}%
\frac{\delta(\xi-\omega_{k,l})+\delta(\xi+\omega_{k,l})}{2}\right]  ,
\notag \end{equation}
where the forward and inverse Fourier transforms of some function $h(t)$ are
defined as follows
\begin{align*}
\hat{h}(\xi)  &  \equiv\left(  \mathcal{F}h\right)  (\xi)\equiv\frac{1}%
{\sqrt{2\pi}}\int\limits_{\mathbb{R}}h(t)e^{-i\xi t}dt,\\
h(t)  &  \equiv\left(  \mathcal{F}^{-1}\hat{h}\right)  (t)\equiv\frac{1}%
{\sqrt{2\pi}}\int\limits_{\mathbb{R}}\hat{h}(\xi)e^{i\xi t}d\xi.
\end{align*}
One can conclude that the Fourier transform $\hat{g}_{2,k}(\xi)$ of
$g_{2,k}(t)$ is the following distribution
\begin{equation}
\hat{g}_{2,k}(\xi)=\sqrt{2\pi}\sum_{l=0}^{\infty}f_{k,l}\frac{\delta
(\xi-\omega_{k,l})+\delta(\xi+\omega_{k,l})}{2} \label{E:deltas}%
\end{equation}
It is tempting to try to recover coefficients $f_{k,l}$ by applying the
Fourier transform $\mathcal{F}$ to each $g_{2,k}(t)$ and thus obtaining
$\hat{g}_{2,k}(\xi)$. However, functions $g_{2,k}(t)$ do not vanish at
infinity, and they (in general) are not periodic on any finite interval
$[0,T]$. Therefore, such a computation would be quite inaccurate, at best.

A technique well known in digital signal processing as a way to Fourier
transform long time series, is to multiply the signal by a window function.
For convenience, let us extend $g_{2,k}(t)$ as an even function to the
interval $[-T,T]$ (formulas (\ref{E:oneside}) and (\ref{E:NUFourier}) will
remain valid under such extension). Consider now an even, infinitely smooth
function $\eta(t)$ vanishing with all its derivatives at $-1$ and $1$, and its
scaled version $\eta_{T}(t)=\eta(t/T)$. Since the product $\eta_{T}%
(t)g_{2,k}(t)$ is periodic on $[-T,T]$, its Fourier transform
\begin{equation}
\widehat{\eta_{T}g_{2,k}}(\xi)\equiv\frac{1}{\sqrt{2\pi}}\int\limits_{-T}%
^{T}\eta_{T}(t)g_{2,k}(t)e^{-i\xi t}dt, \label{E:productFourier}%
\end{equation}
can be easily computed (for example, by applying the composite trapezoid rule
to the discretized signal). Now, by the convolution theorem
\begin{equation}
\widehat{\eta_{T}g_{2,k}}(\xi)=\frac{1}{\sqrt{2\pi}}\widehat{\eta_{T}}%
(\xi)\ast\widehat{g_{2,k}}(\xi)=\frac{1}{\sqrt{2\pi}}\int\limits_{\mathbb{R}%
}\widehat{\eta_{T}}(\zeta-\xi)\widehat{g_{2,k}}(\zeta)d\zeta,
\notag \end{equation}
and, taking into account the equality $\widehat{\eta_{T}}(\xi)=T\widehat{\eta
}(T\xi)$ and equation (\ref{E:deltas}) one obtains
\begin{equation}
\widehat{\eta_{T}g_{2,k}}(\xi)=\sum_{m=0}^{\infty}f_{k,m}\frac{\widehat{\eta
_{T}}(\xi-\omega_{k,m})+\widehat{\eta_{T}}(\xi+\omega_{k,m})}{2}=T\sum
_{m=0}^{\infty}f_{k,m}\frac{\widehat{\eta}(T(\xi-\omega_{k,m}))+\widehat{\eta
}(T(\xi+\omega_{k,m}))}{2}%
\notag \end{equation}

Function $\widehat{\eta_{T}}(\xi)$ is infinitely smooth (as a Fourier
transform of a finitely supported function), and vanishes at infinity faster
than any rational function in $\xi$ (since $\eta(t)$ is infinitely smooth).
One may view the problem of reconstruction of $\hat{g}_{2,k}(\xi)$ from
$\widehat{\eta_{T}}(\xi)\ast\widehat{g_{2,k}}(\xi)$ as a deconvolution
problem. Due to the smoothness of the convolution kernel $\widehat{\eta_{T}%
}(\xi)$ such a deconvolution is an extremely ill-posed problem.

Fortunately, the problem at hand can be solved without resorting to standard
deconvolution techniques. Since we know a priori that the $\hat{g}_{2,k}(\xi)$
is a combination of delta-functions (with unknown coefficients $f_{k,l})$
whose positions in the frequency space are known (see
equation\ (\ref{E:deltas})), the problem becomes much simpler. Let us consider
the values of $\widehat{\eta_{T}g_{2,k}}$ only at the frequencies
$\omega_{k,l}$. Then, for each fixed $k$ we obtain the following system of
linear equations with respect to unknown $f_{k,m}$:
\begin{equation}
T\sum_{m=0}^{\infty}f_{k,m}\frac{\widehat{\eta}(T(\omega_{k,l}-\omega
_{k,m}))+\widehat{\eta}(T(\omega_{k,l}+\omega_{k,m}))}{2}=\widehat{\eta
_{T}g_{2,k}}(\omega_{k,l}),\qquad l=0,1,2,... \label{E:system}%
\end{equation}
Further, by separating the `diagonal' terms this system can be re-written in
the form
\begin{equation}
f_{k,l}+\frac{\widehat{\eta}(2T\omega_{k,l})}{\widehat{\eta}(0)}f_{k,l}%
+\sum_{\underset{m\neq l}{m=0}}^{\infty}\frac{\widehat{\eta}(T(\omega
_{k,l}-\omega_{k,m}))+\widehat{\eta}(T(\omega_{k,l}+\omega_{k,m}%
))}{\widehat{\eta}(0)}f_{k,m}=\frac{2}{T}\frac{\widehat{\eta_{T}g_{2,k}%
}(\omega_{k,l})}{\widehat{\eta}(0)} \label{E:system1}%
\end{equation}
for $l=0,1,2,...,$ $k=0,1,2...,$ except the case when $k=l=0$ which assumes a
simple form%
\begin{equation}
f_{0,0}+\sum_{\underset{m\neq l}{m=0}}^{\infty}\frac{\widehat{\eta}%
(T\omega_{0,m})}{\widehat{\eta}(0)}f_{0,m}=\frac{1}{T}\frac{\widehat{\eta
_{T}g_{2,0}}(\omega_{0,0})}{\widehat{\eta}(0)}. \label{E:system1a}%
\end{equation}

Let us discuss the properties of these systems (one system for every value of
$k)$. Due to the fast decrease of $\widehat{\eta_{T}}(\xi)$ at infinity, as
$T$ goes to infinity, values of $\widehat{\eta}(T(\omega_{k,l}-\omega_{k,m}))$
and $\widehat{\eta}(T(\omega_{k,l}+\omega_{k,m}))$ converge to zero, and the
system becomes diagonal. This suggests that for sufficiently large values of
$T$ one can try to approximately solve the problem by the formula%
\begin{align}
f_{k,l}  &  \thickapprox\frac{2}{T}\frac{\widehat{\eta_{T}g_{2,k}}%
(\omega_{k,l})}{\widehat{\eta}(0)+\widehat{\eta}(2T\omega_{k,l})},\qquad
k,l=0,1,2,...,\label{E:approx-diag}\\
f_{0,0}  &  \thickapprox\frac{1}{T}\frac{\widehat{\eta_{T}g_{2,0}}%
(\omega_{0,0})}{\widehat{\eta}(0)}.
\end{align}
and thus to reconstruct $p_{0}(x)$ from the measurements made on one side of
the square~$\Omega$.

Examples of such approximate reconstructions can be found in~\cite{CHK};
it is clear that the images are distorted by noticeable artifacts. A closer look at
the functions $\widehat{\eta}(T(\omega_{k,l}-\omega_{k,m}))$ reveals the
reason for these distortions. The convergence of these functions to zero
depends on the values of the difference $\omega_{k,l}-\omega_{k,m}$ which is
not uniform with respect to $k$ and $l$. In particular, for large values of
$k$ and $l=0$ and $m=1$ this difference can be approximately estimated as
follows
\begin{equation}
\omega_{k,1}-\omega_{k,0}=\pi\left(  \sqrt{k^{2}+1}-\sqrt{k^{2}}\right)  =\pi
k\left(  \sqrt{1+1/k^{2}}-1\right)  \thickapprox\frac{\pi}{2k},
\notag \end{equation}
so that the value of $T(\omega_{k,l}-\omega_{k,m})$ is not large for large
values of $k$. This implies that for any (large) value of $T$ there are some
values of $k,l$ and $m$ for which off-diagonal terms are not small, and
approximation (\ref{E:approx-diag}) is not accurate. Moreover, one may suspect
(and numerical simulations confirm this) that the equations with large values
of $k$ and small values of $l$ are a source of instability, so that attempts
to solve system (\ref{E:system1}) numerically (for large $k)$ lead to a
significant amplification of the noise present in the data.

\subsection{Reconstruction using from data measured on two adjacent sides}

As shown below, one can eliminate this instability by using not all equations
(\ref{E:system1}) but only those for $l=k,k+1,k+2,...$ . The missing
information can be obtained from the data $g_{1}(x_{2},t)$ measured on the
second side of the square corresponding to $x_{1}=0$. Define functions
$g_{1,l}(t)$ by the equations%
\begin{equation}
g_{1,l}(t)\equiv2\int\limits_{0}^{1}g_{1}(x_{2},t)\cos(\pi lx_{2})dx_{2}.
\notag \end{equation}
The derivation similar to the one in the previous section yields the following
system of equations:%
\begin{equation}
f_{k,l}+\widehat{\eta}(2T\omega_{k,l})f_{k,l}+\sum_{\underset{m\neq k}{m=0}%
}^{\infty}\frac{\widehat{\eta}(T(\omega_{k,l}-\omega_{m,l}))+\widehat{\eta
}(T(\omega_{k,l}+\omega_{m,l}))}{\widehat{\eta}(0)}f_{k,m}=\frac{2}{T}%
\frac{\widehat{\eta_{T}g_{1,l}}(\omega_{k,l})}{\widehat{\eta}(0)},
\label{E:system2}%
\end{equation}
$k,l=0,1,2,...,$ (with a simpler formula for the case $k=l=0$ which we will
omit). We, however, are going to use only half of these equations, namely
those corresponding to $k=l+1,l+2,l+3,...$ in combination with a half of the
equations (\ref{E:system1}), arriving at the following set of equations:%
\begin{equation}
\left\{
\begin{array}
[c]{c}%
f_{0,0}+\sum_{\underset{m\neq l}{m=0}}^{\infty}\frac{\widehat{\eta}%
(T(\omega_{0,m}))}{\widehat{\eta}(0)}f_{0,m}=\frac{1}{T}\frac{\widehat{\eta
_{T}g_{2,0}}(\omega_{0,0})}{\widehat{\eta}(0)},\quad l=k=0,\\
f_{k,l}+\widehat{\eta}(2T\omega_{k,l})f_{k,l}+\sum_{\underset{m\neq l}{m=0}%
}^{\infty}\frac{\widehat{\eta}(T(\omega_{k,l}-\omega_{k,m}))+\widehat{\eta
}(T(\omega_{k,l}+\omega_{k,m}))}{\widehat{\eta}(0)}f_{k,m}=\frac{2}{T}%
\frac{\widehat{\eta_{T}g_{2,k}}(\omega_{k,l})}{\widehat{\eta}(0)},\quad l\geq
k,\\
f_{k,l}+\widehat{\eta}(2T\omega_{k,l})f_{k,l}+\sum_{\underset{m\neq k}{m=0}%
}^{\infty}\frac{\widehat{\eta}(T(\omega_{k,l}-\omega_{m,l}))+\widehat{\eta
}(T(\omega_{k,l}+\omega_{m,l}))}{\widehat{\eta}(0)}f_{k,m}=\frac{2}{T}%
\frac{\widehat{\eta_{T}g_{1,l}}(\omega_{k,l})}{\widehat{\eta}(0)},\quad k>l.
\end{array}
\right.  \label{E:big-system}%
\end{equation}
Since for very large values of $T$ factors $\widehat{\eta}(T(...))$ vanish
(while $\widehat{\eta}(0)$ is a constant), one can compute a crude
approximation $f^{(0)}(x)$ by defining the coefficients $F^{(0)}$ as follows%
\begin{align}
f_{k,l}^{(0)}  &  =\left\{
\begin{array}
[c]{cc}%
\frac{1}{T}\frac{\widehat{\eta_{T}g_{2,0}}(\omega_{0,0})}{\widehat{\eta}%
(0)}, & k=l=0,\\
\frac{2}{T}\frac{\widehat{\eta_{T}g_{2,k}}(\omega_{k,l})}{\widehat{\eta}%
(0)}, & l\geq k\\
\frac{2}{T}\frac{\widehat{\eta_{T}g_{1,l}}(\omega_{k,l})}{\widehat{\eta}%
(0)}, & k>l
\end{array}
\right.  ,\label{E:first-appro1}\\
F^{(0)}  &  =\left\{  f_{k,l}^{(0)}\right\}  _{k,l=0}^{\infty},
\label{E:first-appro2}%
\end{align}
and by summing the Fourier series%
\begin{equation}
f^{(0)}(x)=\left[  \mathcal{F}_{\operatorname{series}}^{-1}F^{(0)}\right]
(x). \label{E:inverseFourier}%
\end{equation}
Formulas (\ref{E:first-appro1}), (\ref{E:first-appro2}) define a linear
operator $\mathcal{R}$ that maps boundary values $\mathbf{g}$ into the set of
Fourier coefficients $F$, so that%
\begin{equation}
f^{(0)}=\mathcal{F}_{\operatorname{series}}^{-1}F^{(0)}=\mathcal{F}%
_{\operatorname{series}}^{-1}\mathcal{R}\mathbf{g=}\mathcal{F}%
_{\operatorname{series}}^{-1}\mathcal{RW}f.
\notag \end{equation}

As our numerical experiments show, approximation $f^{(0)}$ yields
qualitatively correct images even for moderate values of $T$. Our goal,
however, is to obtain a quantitatively accurate, theoretically exact
reconstruction. Forming and decomposing a matrix corresponding to the system
(\ref{E:big-system}) (after proper truncation of the spatial harmonics) is
computationally prohibitive, since the number of unknowns in the 3D high
resolution images can reach hundreds of millions. Instead, we will solve this
system iteratively.

Let us define an operator $\mathcal{A}$ that transforms a two-dimensional
sequence of numbers $H=\{h_{k,m}\}_{k,m=0}^{\infty}$ into a two-dimensional
sequence of numbers $E=\{e_{k,m}\}_{k,m=0}^{\infty}$ by the formulas%
\begin{equation}%
\begin{array}
[c]{c}%
e_{0,0}=\sum_{\underset{m\neq l}{m=0}}^{\infty}\frac{\widehat{\eta}%
(T(\omega_{0,m}))}{\widehat{\eta}(0)}h_{0,m},\quad l=k=0,\\
e_{k,l}=\frac{\widehat{\eta}(2T\omega_{k,l})}{\widehat{\eta}(0)}h_{k,l}%
+\sum_{\underset{m\neq l}{m=0}}^{\infty}\frac{\widehat{\eta}(T(\omega
_{k,l}-\omega_{k,m}))+\widehat{\eta}(T(\omega_{k,l}+\omega_{k,m}%
))}{\widehat{\eta}(0)}h_{k,m},\quad l\geq k,\\
e_{l,k}=\frac{\widehat{\eta}(2T\omega_{k,l})}{\widehat{\eta}(0)}h_{k,l}%
+\sum_{\underset{m\neq k}{m=0}}^{\infty}\frac{\widehat{\eta}(T(\omega
_{k,l}-\omega_{m,l}))+\widehat{\eta}(T(\omega_{k,l}+\omega_{m,l}%
))}{\widehat{\eta}(0)}h_{k,m},\quad k>l.
\end{array}
\label{E:big-oper}%
\end{equation}
or, in the operator notation,%
\begin{equation}
E=\mathcal{A}H.
\notag \end{equation}
Then equations (\ref{E:big-system}) can be re-written as%
\begin{equation}
(\mathcal{I+A)}F=F^{(0)}, \label{E:coef-eqs}%
\end{equation}
where $F$ is the set of the Fourier coefficients of the sought $f(x),$ and
$F^{(0)}$ are the Fourier coefficients of the initial, crude approximation
$f^{(0)}(x)$. If the spectral radius of operator $\mathcal{A}$ is less than 1,
the unique solution of\ equation (\ref{E:coef-eqs}) can be found in the form
of a converging Neumann series
\begin{equation}
F=\sum\limits_{K=0}^{\infty}(-\mathcal{A})^{K}F^{(0)}. \label{E:Neumann-coef}%
\end{equation}
Function $f(x)$ is then found from $F$ by summing the Fourier series (see
(\ref{E:inverseFourier})). Alternatively, this solution ($F$) can be
represented as the limit of iterations defined by the recurrence relation%
\begin{equation}
F^{(K)}=F^{(0)}-\mathcal{A}F^{(K-1)},\quad K=1,2,3,... \label{E:Fourier-iter}%
\end{equation}
In section~\ref{S:convergence} we analyze operator $\mathcal{A}$ and find a sufficient
condition for the convergence of the Neumann series (\ref{E:Neumann-coef}).

In practical computations we are not evaluating series (\ref{E:Neumann-coef})
directly, since this is still expensive. By summing the Fourier series (i.e.
by applying the $\mathcal{F}_{\operatorname{series}}^{-1}$ operator)
iterations (\ref{E:Fourier-iter}) can be replaced by the equivalent relation%
\begin{equation}
f^{(K)}(x)=f^{(0)}(x)-\mathcal{F}_{\operatorname{series}}^{-1}\mathcal{AF}%
_{\operatorname{series}}f^{(K-1)},\quad K=1,2,3,...
\notag \end{equation}
or%
\begin{equation}
f^{(K)}(x)=\mathcal{F}_{\operatorname{series}}^{-1}\mathcal{R}\mathbf{g}%
-\mathcal{F}_{\operatorname{series}}^{-1}\mathcal{AF}_{\operatorname{series}%
}f^{(K-1)},\quad K=1,2,3,... \label{E:space-iter}%
\end{equation}
On the other hand, it follows from (\ref{E:coef-eqs}) that%
\begin{equation}
(\mathcal{I+A)F}_{\operatorname{series}}f=\mathcal{RW}f,
\notag \end{equation}
and, since (\ref{E:coef-eqs}) holds for any initial condition $f(x),$ the
following operator identity holds%
\begin{equation}
\mathcal{F}_{\operatorname{series}}^{-1}(\mathcal{I+A)F}%
_{\operatorname{series}}=\mathcal{F}_{\operatorname{series}}^{-1}%
\mathcal{RW}.
\notag \end{equation}
Thus,%
\begin{equation}
-\mathcal{F}_{\operatorname{series}}^{-1}\mathcal{AF}_{\operatorname{series}%
}=\mathcal{I-F}_{\operatorname{series}}^{-1}\mathcal{RW},
\notag \end{equation}
and, therefore, recurrence relation (\ref{E:space-iter}) is equivalent to
\begin{equation}
f^{(K)}=\mathcal{F}_{\operatorname{series}}^{-1}\mathcal{R}\mathbf{g+}%
f^{(K-1)}-\mathcal{F}_{\operatorname{series}}^{-1}\mathcal{RW}f^{(K-1)},\quad
K=1,2,3,...
\notag \end{equation}
or to%
\begin{align}
f^{(K)}  &  =f^{(K-1)}+\mathcal{F}_{\operatorname{series}}^{-1}\mathcal{R(}%
\mathbf{g}-\mathcal{W}f^{(K-1)}),\quad K=1,2,3,...,\label{E:space-iter1}\\
f^{(0)}  &  =\mathcal{F}_{\operatorname{series}}^{-1}\mathcal{R}\mathbf{g.}
\label{E:space-iter2}%
\end{align}
Our computational algorithm is based on equations (\ref{E:space-iter1}),
(\ref{E:space-iter2}); these iterations are theoretically equivalent to those
defined by equation (\ref{E:Fourier-iter}) and have the same convergence
properties. The computational advantage in using (\ref{E:space-iter1}),
(\ref{E:space-iter2})
lies in the possibility of easily implementing all the necessary
operations using fast transforms (FFTs). The implementation
details are discussed in Section~\ref{S:numerics}.
In the next section we analyze the convergence of such iterations,
and the form (\ref{E:Fourier-iter}) is more convenient for this purpose.

\subsection{Convergence analysis\label{S:convergence}}

Let us consider the $||\cdot||_{\infty}$ norm defined on the space of
double-indexed infinite sequences $H=\{h_{k,l}\}_{k,l=0}^{\infty}$:%
\begin{equation}
||H||_{\infty}=\sup_{k,l=0.1,2,3...}|h_{k,l}|,
\notag \end{equation}
and a space $\mathbb{L}$ of such sequences with bounded $||\cdot||_{\infty}$
norm. Further, for a linear operator $\mathcal{B}:\mathbb{L\rightarrow L}$
define the induced $\infty$- norm:%
\begin{equation}
||\mathcal{B}||_{\infty}=\sup_{||H||_{\infty}\neq0}\frac{||\mathcal{B}%
(H)||_{\infty}}{||H||_{\infty}}.
\notag \end{equation}
Our goal is to estimate the induced $\infty$- norm of the operator
$\mathcal{A}$ defined by equations (\ref{E:big-oper}) in the previous
sections.\ These equations contain factors in the form $\widehat{\eta
}(T(\omega_{k,l}-\omega_{k,m}))$. Note that for a sufficiently smooth cut-off
function $\eta(t)$ its Fourier transform declines fast at infinity:
\begin{equation}
|\widehat{\eta}(\xi)|\leq\frac{B}{1+|\xi|^{2}}.
\notag \end{equation}
Let us find a lower bound on $|\omega_{k,l}-\omega_{k,m}|$. Consider first the
case $l\geq k$:%
\begin{equation}
|\omega_{k,l}-\omega_{k,m}|=\pi\left\vert \sqrt{k^{2}+l^{2}}-\sqrt{k^{2}%
+m^{2}}\right\vert =\frac{\pi|l-m|(l+m)}{\sqrt{k^{2}+l^{2}}+\sqrt{k^{2}+m^{2}%
}}.
\notag \end{equation}
If, in addition $m\geq k,$ the above equation yields%
\begin{equation}
|\omega_{k,l}-\omega_{k,m}|\geq\pi\frac{|l-m|(l+m)}{\sqrt{2}(l+m)}=\frac
{\pi|l-m|}{\sqrt{2}}.
\notag \end{equation}
If $m<k$ (and $l\geq k)$:%
\begin{equation}
|\omega_{k,l}-\omega_{k,m}|\geq\frac{\pi|l-m|(l+m)}{\sqrt{2}(l+k)}\geq
\frac{\pi|l-m|(l+m)}{2\sqrt{2}l}\geq\frac{\pi|l-m|}{2\sqrt{2}},
\notag \end{equation}
so that the uniform bound holds if $l\geq k$ for all values of $m$:%
\begin{equation}
|\omega_{k,l}-\omega_{k,m}|\geq\frac{\pi|l-m|}{2\sqrt{2}}.
\notag \end{equation}
Now one can bound $\left\vert \widehat{\eta}(T(\omega_{k,l}-\omega
_{k,m}))\right\vert $, still under assumption $l\geq k$:%
\begin{align}
\left\vert \widehat{\eta}(T(\omega_{k,l}-\omega_{k,m}))\right\vert  &
\leq\frac{B}{1+T^{2}(\omega_{k,l}-\omega_{k,m})^{2}}\leq\frac{B}{1+\pi
^{2}T^{2}(l-m)^{2}/8}\nonumber\\
&  \leq\frac{B}{\pi^{2}T^{2}(l-m)^{2}/8},\text{ if }m\neq l. \label{E:est1}%
\end{align}
Similarly, if $l>0$ or $k>0,$%
\begin{equation}
\left\vert \widehat{\eta}(T(\omega_{k,l}+\omega_{k,m}))\right\vert \leq
\frac{B}{1+T^{2}(\omega_{k,l}+\omega_{k,m})^{2}}\leq\frac{B}{1+T^{2}\pi
^{2}(1+m^{2})}\leq\frac{B}{T^{2}\pi^{2}(1+m^{2})} \label{E:est2}%
\end{equation}
and
\begin{equation}
\left\vert \widehat{\eta}(2T\omega_{k,l})\right\vert \leq\frac{B}{T^{2}\pi
^{2}}. \label{E:est3}%
\end{equation}

Let us now estimate the $\infty$-norm of the vector $G=\mathcal{A}H$ resulting
from applying operator $\mathcal{A}$ to a sequence $H$ whose $\infty$-norm
equals 1, so that%
\begin{equation}
|h_{k,l}|\leq1,\qquad l,k=0,1,2,...\text{ .} \label{E:est4}%
\end{equation}
First, we use the second equation in (\ref{E:big-oper}) and bound $e_{k,l}$
with $l\geq k$ (excluding the case $k=l=0)$. Taking into account inequalities
(\ref{E:est1})-(\ref{E:est4}) one obtains%
\begin{align*}
|e_{k,l}|  &  \leq\frac{B}{T^{2}\pi^{2}}\left(  1+\sum_{\underset{m\neq
l}{m=0}}^{\infty}\frac{1}{(1+m^{2})}+\frac{8}{(l-m)^{2}}\right) \\
&  \leq\frac{B}{T^{2}\pi^{2}}\left(  1+\sum_{m=0}^{\infty}\frac{1}{(1+m^{2}%
)}+16\sum_{m=1}^{\infty}\frac{1}{m^{2}}\right) \\
&  \leq\frac{B}{T^{2}\pi^{2}}\left(  2+17\sum_{m=1}^{\infty}\frac{1}{m^{2}%
}\right)  .
\end{align*}
The value of the latter series is well known (see~\cite{Gr-Ryzhik}, equation 0.233)%
\begin{equation}
\sum_{m=1}^{\infty}\frac{1}{m^{2}}=\frac{\pi}{6},
\notag \end{equation}
resulting in the estimate%
\begin{equation}
|e_{k,l}|\leq\frac{B}{T^{2}\pi^{2}}\frac{12+17\pi^{2}}{6}. \label{E:nice-est}%
\end{equation}
The above computations can be replicated with very minor changes to show that
(\ref{E:nice-est}) holds also in the cases $k=l=0$ and $k>l$ (corresponding to
the first and third lines of (\ref{E:big-oper})) which proves the following

\begin{lemma}
Operator $\mathcal{A}$ is bounded in the induced $\infty$-norm by%
\begin{equation}
||\mathcal{A}||_{\infty}\leq\frac{B}{T^{2}\pi^{2}}\frac{12+17\pi^{2}}{6}.
\notag \end{equation}

\end{lemma}

The Lemma implies that if the acquisition time $T$ is large enough, e.g. if
\begin{equation}
T>\frac{1}{\pi}\sqrt{\frac{B(12+17\pi^{2})}{6}}, \label{E:bigtime}%
\end{equation}
then operator $\mathcal{A}$ is a contraction mapping in the induced $\infty
$-norm, i.e. $||\mathcal{A}||_{\infty}<1$. In this case the standard theory of
contraction mappings applies and one arrives at the following

\begin{theorem}
For a sufficiently large acquisition time $T$ (satisfying (\ref{E:bigtime}))
the Neumann series (\ref{E:Neumann-coef}) converges in the $\infty$-norm,
implying convergence of iterations (\ref{E:space-iter1}), (\ref{E:space-iter2}%
) in the following sense%
\begin{equation}
\sup_{k,l=0.1,2,3...}|f_{k,l}^{(K)}-f_{k,l}|\underset{K\rightarrow
\infty}{\rightarrow}0.
\notag \end{equation}

\end{theorem}

\bigskip

\subsection{3D version of the method}

In the 3D case the problem and the proposed algorithm are very similar to
their 2D counterparts. Since the 3D case is both more important for the
practical use and more challenging computationally, we provide in this section
a brief outline of the derivation and convergence analysis of our technique.

The acoustic pressure $p(t,x),$ $x=(x_{1},x_{2},x_{3})$ solves IBVP
(\ref{E:Forward1})-(\ref{E:Forward3}) in 3D, in a cube $\Omega=[0,1]\times
\lbrack0,1]\times\lbrack0,1]$. As before, the speed of sound $c_{0}$ will be
assumed to equal 1. The measurements are represented by the triple of time
series $\mathbf{g}=(g_{1}(t,x_{2}),g_{2}(t,x),g_{3}(t,x))$ corresponding to
the pressure values on the three adjacent sides of the cube $\Omega$:%
\begin{equation}
g_{1}(t,x_{2},x_{3})=\left.  p(t,x)\right\vert _{x_{1}=0},\ g_{2}%
(t,x_{1},x_{3})=\left.  p(t,x)\right\vert _{x_{2}=0},\ g_{3}(t,x_{1}%
,x_{2})=\left.  p(t,x)\right\vert _{x_{3}=0},\ t\in\lbrack0,T].
\notag \end{equation}
The solution of IBVP in the cube can be represented using the Neumann
Laplacian eigenfunctions%
\begin{equation}
\varphi_{k,l,n}(x)=\cos(\pi kx_{1})\cos(\pi lx_{2})\cos(\pi nx_{3})
\notag \end{equation}
with eigenfrequencies
\begin{equation}
\omega_{k,l,n}=\pi\sqrt{k^{2}+l^{2}+n^{2}},\quad k,l,n=0,1,2,...\text{ .}%
\notag \end{equation}
It has the following form:%
\begin{equation}
p(t,x)=\sum_{k=0}^{\infty}\sum_{l=0}^{\infty}\sum_{n=0}^{\infty}f_{k,l,n}%
\cos(\pi kx_{1})\cos(\pi lx_{2})\cos(\pi nx_{3})\cos(\omega_{k,l,n}t),\quad
x\in\Omega,\quad t\in\lbrack0,T].
\notag \end{equation}
Coefficients $F=\{f_{k,l,n}\}_{k,l,n=0}^{\infty}$ are related to $f(x)$
through the Fourier cosine series, \ $F=\mathcal{F}_{\operatorname{series}}f$.
They are computed as follows:%
\begin{equation}
f_{k,l,n}=\frac{1}{||\varphi_{k,l,n}||_{2}^{2}}\int\limits_{\Omega}%
f(x)\varphi_{k,l,n}(x)~dx,\quad k,l,n=0,1,2,...\text{ .} \label{E:Fourier3D}%
\end{equation}
Conversely, $f$ is obtained from $F$ by summing the standard 3D cosine Fourier
series, $f=\mathcal{F}_{\operatorname{series}}^{-1}F$.

Consider side of the cube corresponding to $x_{2}=0$. Then%
\begin{equation}
g_{2}(t,x_{1},x_{3})\equiv\left.  p(t,x)\right\vert _{x_{2}=0}=\sum
_{k=0}^{\infty}\sum_{n=0}^{\infty}\cos(\pi kx_{1})\cos(\pi nx_{3})\left[
\sum_{l=0}^{\infty}f_{k,l,n}\cos(\omega_{k,l,n}t)\right],\ t\in
\lbrack0,T]. \label{E:oneside3D}%
\end{equation}
Due to the orthogonality of the products of cosines, the functions in brackets
(denote them $g_{2,k,n}$) can be easily found from $g_{2}$ by the 2D cosine
Fourier series.%
\begin{equation}
g_{2,k,n}(t)\equiv\sum_{l=0}^{\infty}f_{k,l,n}\cos(\omega_{k,l,n}t),\quad
t\in\lbrack0,T],\quad k,n=0,1,2,...\text{ .} \label{E:onesideNUFFT}%
\end{equation}
By finding the values of the Fourier transform $\widehat{g_{2,k,n}\eta_{T}}$
of the product $g_{2,k,n}(t)\eta_{T}(t)$ at the frequencies $\omega_{k,l,n}$
one obtains the following infinite system of equations%
\begin{equation}
f_{k,l,n}+\frac{\widehat{\eta}(2T\omega_{k,l,n})}{\widehat{\eta}(0)}%
f_{k,l,n}+\sum_{\underset{m\neq l}{m=0}}^{\infty}\frac{\widehat{\eta}%
(T(\omega_{k,l,n}-\omega_{k,m,n}))+\widehat{\eta}(T(\omega_{k,l,n}%
+\omega_{k,m,n}))}{\widehat{\eta}(0)}f_{k,m,n}=\frac{2}{T}\frac{\widehat{\eta
_{T}g_{2,k,n}}(\omega_{k,l,n})}{\widehat{\eta}(0)},
\notag \end{equation}
for $k,l,n=0,1,2,...,$ except the case when $k=l=n=0$ which yields a simpler
equation%
\begin{equation}
f_{0,0,0}+\sum_{\underset{m\neq l}{m=0}}^{\infty}\frac{\widehat{\eta}%
(T\omega_{0,m,0})}{\widehat{\eta}(0)}f_{0,m,0}=\frac{1}{T}\frac{\widehat{\eta
_{T}g_{2,0,0}}(\omega_{0,0})}{\widehat{\eta}(0)}.
\notag \end{equation}
As in the 2D case, for a stable reconstruction we need to use equations from
three sides (equations for the remaining two sides are obtain in a similar
fashion). We thus arrive at the system%
\begin{equation}
\left\{
\begin{array}
[c]{c}%
f_{0,0,0}+\left[  \sum_{\underset{m\neq l}{m=0}}^{\infty}\frac{\widehat{\eta
}(T\omega_{0,m,0})}{\widehat{\eta}(0)}f_{0,m,0}\right]  =f_{0,0,0}^{(0)},\quad
l=k=n=0,\\
f_{k,l,n}+\left[  \frac{\widehat{\eta}(2T\omega_{k,l,n})}{\widehat{\eta}%
(0)}f_{k,l,n}+\sum_{\underset{m\neq l}{m=0}}^{\infty}\frac{\widehat{\eta
}(T(\omega_{k,l,n}-\omega_{m,l,n}))+\widehat{\eta}(T(\omega_{k,l,n}%
+\omega_{m,l,n}))}{\widehat{\eta}(0)}f_{m,l,n}\right]  =f_{k,l,n}%
^{(0)},\,k\geq l,\,k\geq n,\\
f_{k,l,n}+\left[  \frac{\widehat{\eta}(2T\omega_{k,l,n})}{\widehat{\eta}%
(0)}f_{k,l,n}+\sum_{\underset{m\neq l}{m=0}}^{\infty}\frac{\widehat{\eta
}(T(\omega_{k,l,n}-\omega_{k,m,n}))+\widehat{\eta}(T(\omega_{k,l,n}%
+\omega_{k,m,n}))}{\widehat{\eta}(0)}f_{k,m,n}\right]  =f_{k,l,n}%
^{(0)},\,l>k,\,l\geq n,\\
f_{k,l,n}+\left[  \frac{\widehat{\eta}(2T\omega_{k,l,n})}{\widehat{\eta}%
(0)}f_{k,l,n}+\sum_{\underset{m\neq l}{m=0}}^{\infty}\frac{\widehat{\eta
}(T(\omega_{k,l,n}-\omega_{k,l,m}))+\widehat{\eta}(T(\omega_{k,l,n}%
+\omega_{k,l,m}))}{\widehat{\eta}(0)}f_{k,l,m}\right]  =f_{k,l,n}%
^{(0)},\,n>l,\,n>k,
\end{array}
\right.
\notag \end{equation}
where the right hand side corresponds to the Fourier coefficients of the crude
first approximation $f^{(0)}(x)$:%
\begin{equation}
f_{k,l,n}^{(0)}=\left\{
\begin{array}
[c]{cc}%
\frac{1}{T}\frac{\widehat{\eta_{T}g_{2,0,0}}(\omega_{0,0})}{\widehat{\eta}%
(0)}, & l=k=n=0,\\
\frac{2\widehat{\eta_{T}g_{1,l,n}}(\omega_{k,l,n})}{T\widehat{\eta}(0)}, &
\,k\geq l,\,k\geq n,\\
\frac{2\widehat{\eta_{T}g_{2,k,n}}(\omega_{k,l,n})}{T\widehat{\eta}(0)}, &
l>k,\,l\geq n,\\
\frac{2\widehat{\eta_{T}g_{3,k,l}}(\omega_{k,l,n})}{T\widehat{\eta}(0)}, &
n>l,\,n>k,
\end{array}
\right.   \label{E:diagonal3D}%
\end{equation}
and where the terms in the brackets define the operator $\mathcal{A}$. Now the
system is in the form (\ref{E:coef-eqs}) and its solution can be found in the
form of the Neumann series (\ref{E:Fourier-iter}) if the series converges.
Under the latter condition the solution to the original inverse problem $f(x)$
can be obtained by the iterative process (\ref{E:space-iter1}),
(\ref{E:space-iter2}), where the operator $\mathcal{R}$ is now defined by
(\ref{E:diagonal3D}) and the operator $\mathcal{W}$ maps solutions of the IBVP
(\ref{E:Forward1})-(\ref{E:Forward3}) in 3D corresponding to the initial
condition $f(x)$ into the values of $p(t,x)$ at the three sides of the cube
$\Omega$.

The convergence analysis is also very similar to that in the 2D case. For example,
the difference $|\omega_{k,l,n}-\omega_{k,m,n}|$ can be estimated as follows:
\begin{equation}
|\omega_{k,l,n}-\omega_{k,m,n}|=\pi\left\vert \sqrt{k^{2}+l^{2}+n^{2}}%
-\sqrt{k^{2}+m^{2}+n^{2}}\right\vert =\pi\frac{|l-m|(l+m)}{\sqrt{k^{2}%
+l^{2}+n^{2}}+\sqrt{k^{2}+m^{2}+n^{2}}}.
\notag \end{equation}
If $l\geq k$ \ and $l\geq n$ and $m\geq l$, then%
\begin{equation}
|\omega_{k,l,n}-\omega_{k,m,n}|\geq\pi\frac{|l-m|(l+m)}{\sqrt{3l^{2}}%
+\sqrt{3m^{2}}}\geq\frac{\pi|l-m|}{\sqrt{3}}.
\notag \end{equation}
If $l\geq k$ \ and $l\geq n$ and $m<l,$ then%
\begin{equation}
|\omega_{k,l,n}-\omega_{k,m,n}|\geq\pi\frac{|l-m|(l+m)}{\sqrt{3l^{2}}%
+\sqrt{3l^{2}}}\geq\pi\frac{|l-m|}{2\sqrt{3}}.
\notag \end{equation}
The latter estimate therefore holds for all values of $m,$ if $l\geq k$ \ and
$l\geq n$. Proceeding the same way as in the 2D case one obtains the following
estimate on the $\infty$-norm of the operator $\mathcal{A}$.

\begin{lemma}
(3D case.) Operator $\mathcal{A}$ is bounded in the induced $\infty$-norm by%
\begin{equation}
||\mathcal{A}||_{\infty}\leq\frac{B}{T^{2}\pi^{2}}\frac{12+25\pi^{2}}{6}.
\notag \end{equation}

\end{lemma}

This immediately leads to a sufficient condition for the convergence of the 3D
algorithm if the acquisition time $T$ satisfies the condition
\begin{equation}
T>\frac{1}{\pi}\sqrt{\frac{B(12+25\pi^{2})}{6}}. \label{E:bigtime3D}%
\end{equation}

\begin{theorem}
(3D case) For a sufficiently large acquisition time $T$ (satisfying
(\ref{E:bigtime3D})) the Neumann series (\ref{E:Neumann-coef}) converges in
the $\infty$-norm, implying convergence of iterations (\ref{E:space-iter1}),
(\ref{E:space-iter2}) in the following sense%
\begin{equation}
\sup_{k,l,n=0.1,2,3...}|f_{k,l,n}^{(K)}-f_{k,l,n}|\underset{K\rightarrow
\infty}{\rightarrow}0.
\notag \end{equation}

\end{theorem}

\bigskip

\section{Numerical implementation and simulations\label{S:numerics}}

In this section the proposed reconstruction algorithm is implemented
numerically and tested in numerical simulations. Since real measurements
are made in 3D, we will concentrate on the 3D version of the method. The 2D
algorithm is very similar.

In 3D tomography reconstructions, numerical complexity becomes a major issue,
since in a fine 3D mesh the number of unknowns can easily exceed hundreds of
millions. Thus, it is highly desirable to have an asymptotically fast
reconstruction algorithm reconstructing images on a $N\times N\times N$
computational grid in (preferably) $\mathcal{O}(N^{3})$ or (at most)
$\mathcal{O}(N^{3}\log N)$ floating point operations (flops). Our technique
achieves the latter operation count by utilizing the various Fast Fourier
Transforms (FFTs) on all computationally intensive steps of the algorithm, as
described below.

\begin{figure}[t]
\begin{center}
\subfigure[]{\includegraphics[width=1.9in,height=1.9in]{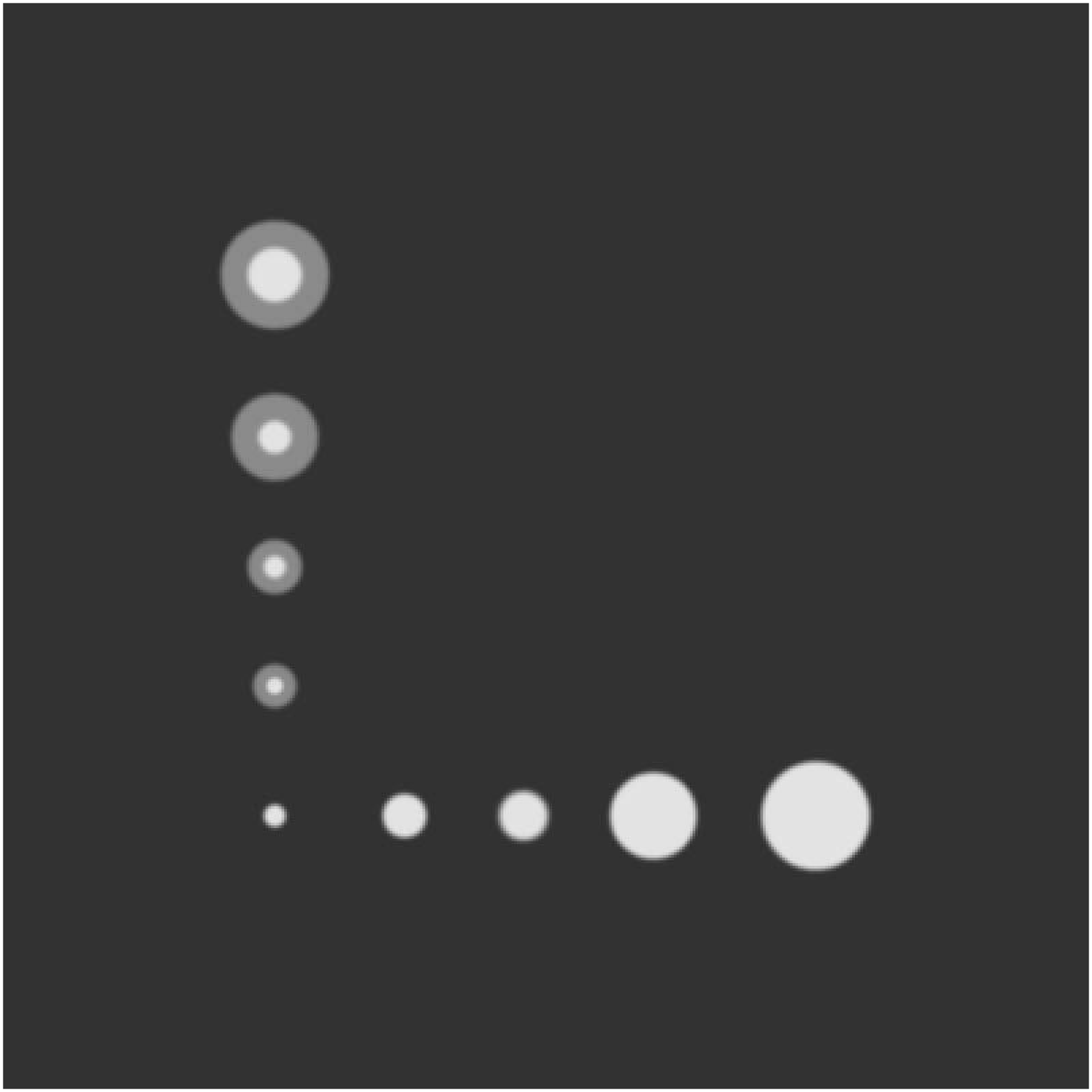}}
\subfigure[]{\includegraphics[width=1.9in,height=1.9in]{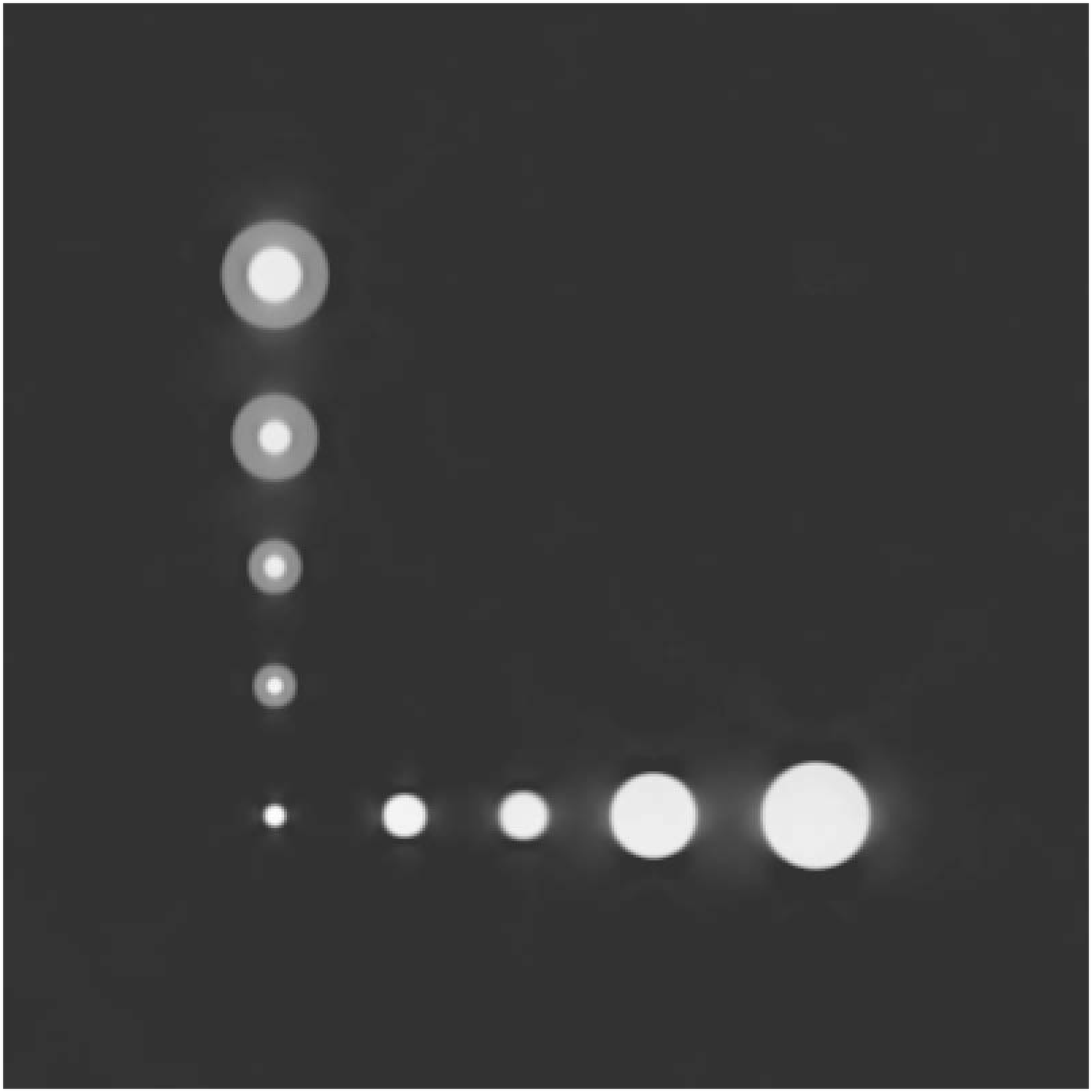}}
\subfigure[]{\includegraphics[width=1.9in,height=1.9in]{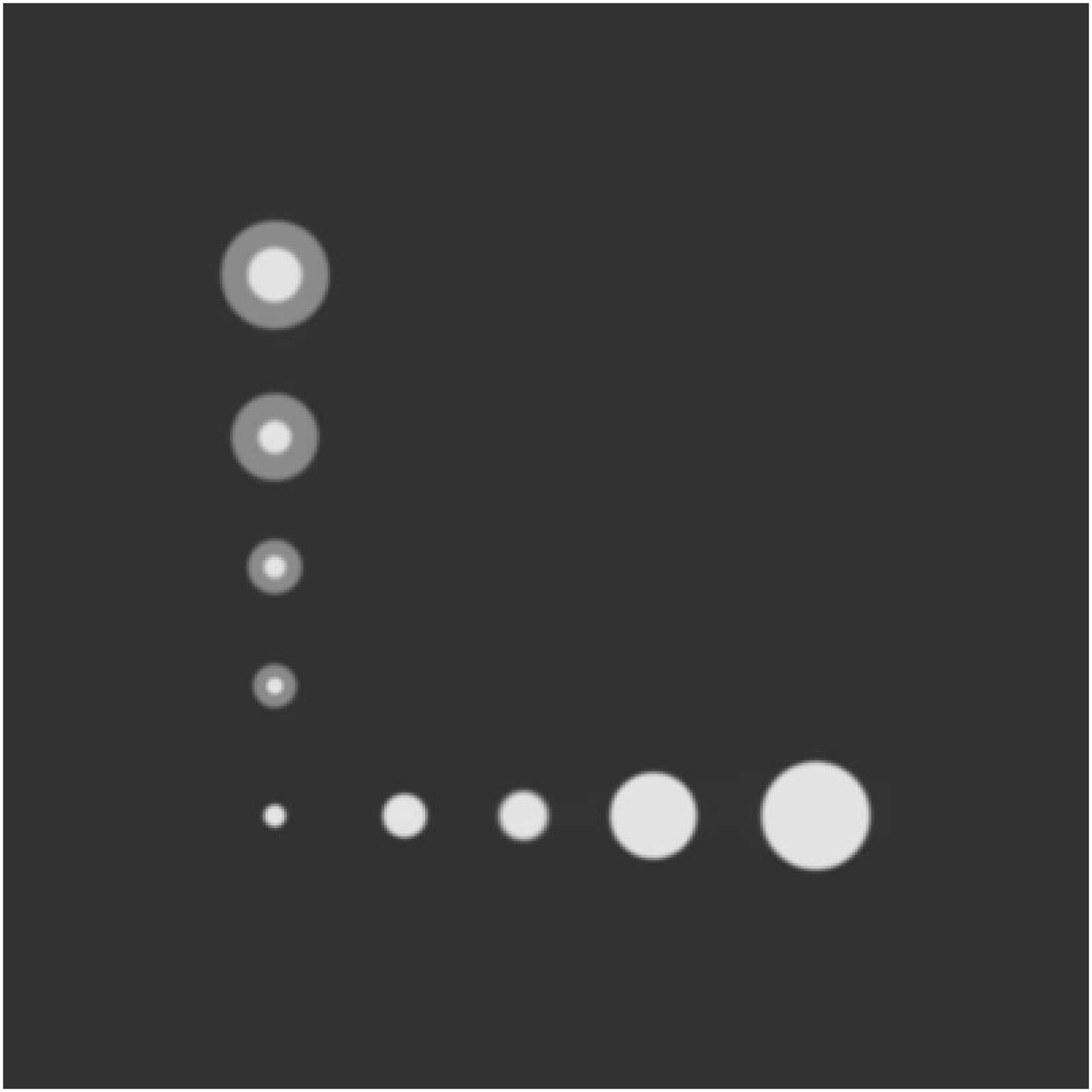}}
\end{center}
\caption{ (a) Phantom, section through the plane $x_3=0.25$ (b) Initial approximation
$f^{(0)}(x)$ (c) 2nd iteration (i.e. $f^{(2)}(x)$) }%
\label{F:n1}%
\end{figure}

\begin{figure}[t]
\begin{center}
\includegraphics[width=5in,height=1.0in]{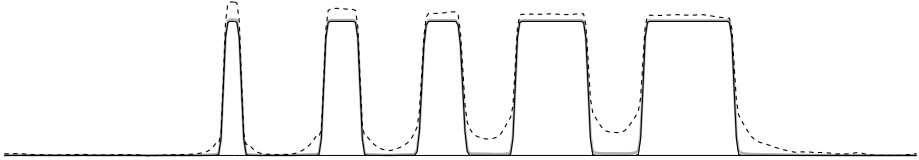}
\end{center}
\caption{A cross-section through the images in Figure~\ref{F:n1} along the line
$x_{2}=0.25$, $x_{3}=0.25$. The solid line represents the phantom, the dashed
line shows $f^{(0)}(x)$ and the gray line represents $f^{(2)}(x)$.}%
\label{F:n2}%
\end{figure}

\subsection{Forward problem}

We remind the reader that the reconstruction algorithm is defined by equations
(\ref{E:space-iter1}), (\ref{E:space-iter2}). $\mathcal{W}$ is one of the
operators in (\ref{E:space-iter1}); it maps the initial condition
$p(0,x)=f(x)$ (here $f(x)$ is an arbitrary function) into values of $p(t,x)$
at the cube faces. In our implementation, computation for each face is done
separately. For example, in order to find $g_{2}(t,x_{1},x_{3})\equiv\left.
p(t,x)\right\vert _{x_{2}=0}$ the following steps are done:

\begin{enumerate}
\item Expand $f(x)$ into the 3D Fourier cosine series, obtain coefficients
$f_{k,l,n}$ (equation (\ref{E:Fourier3D})), $k,l,n=0,1,2,...N$.

\item On a uniform grid in $t,$ compute functions $g_{2,k,n}(t),$
$k,n=0,1,2,...N,$ (equation (\ref{E:onesideNUFFT})).

\item For each value of $t,$ find $g_{2}(t,x_{1},x_{3})$ by summing 2D cosine
Fourier series (\ref{E:oneside3D}).
\end{enumerate}

Let us obtain a crude estimate of the number of operations involved in such a
computation. For simplicity, in all simulations we kept the grid step in time
and in space the same. Therefore, for moderate measurement times $T$ the
number of time nodes is $\mathcal{O}(N)$. Step 1 is done by the 3D cosine FFT,
it requires $\mathcal{O}(N^{3}\log N)$ flops.

Step 2 is slightly more complicated. Evaluation of (\ref{E:onesideNUFFT}) can
be viewed as computing on the uniform grid the 1D Fourier transform of a
function given on a non-uniform grid (frequencies $\omega_{k,l,n}$ are not
equispaced!). This can be done fast by applying the 1D Nonequispaced Fast
Fourier Transform (NSFFT, see, for example \cite{Rokhlin}). This algorithm is
not exact (unlike the regular FFT). However, any needed accuracy can be
obtained, at the expense of increased constant factor implicit in the
$\mathcal{O}(N\log N)$ estimate of the complexity of this technique. In our
simulations the accuracy of NSFFT was of order of $10^{-5}$. Since there is
$N\times N$ of such transforms, the complexity of Step 2 is $\mathcal{O}%
(N^{3}\log N)$ flops.

Step 3 is implemented by summing the 2D cosine series (using the corresponding
type of FFT of size $N\times N$, $\mathcal{O}(N^{2}\log N)$ flops each) for
each of $\mathcal{O}(N)$ time nodes. This step is clearly done in
$\mathcal{O}(N^{3}\log N)$ flops.

To summarize the expense of finding $g_{2}(t,x_{1},x_{3})$ by our method is
$\mathcal{O}(N^{3}\log N)$ flops. The data for the other faces of the cube are
computed in a similar fashion (Step 1 is exactly the same and is performed
only once).

We also used this algorithm to compute the simulated measurements $\mathbf{g}$
for our experiments. In some of the simulations a significant level of noise was added
to $\mathbf{g}$ to model the measurement errors and to check the
stability of the algorithm to such errors.
\noindent\begin{figure}[t]
\begin{center}
\includegraphics[width=5in,height=1.0in]{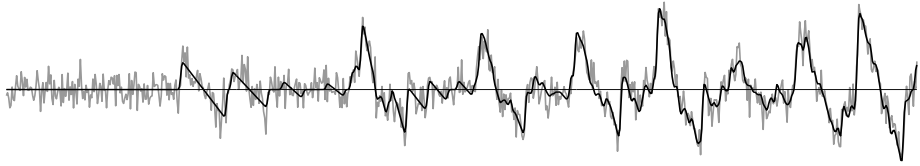}
\end{center}
\caption{Illustration of $100\%$ noise: black line shows exact signal
$g_{1}(t,0.5,0.5)$. Gray line represents the same signal with added noise}%
\label{F:n3}%
\end{figure}

\subsection{Approximate inversion}
The iterations given by equation (\ref{E:space-iter1}) can be viewed as
repeated applications of the forward operator $\mathcal{W}$ and approximate inverse, $\mathcal{F}%
_{\operatorname{series}}^{-1}\mathcal{R}$, in alternating order. The summation
of the Fourier series $\mathcal{F}_{\operatorname{series}}^{-1}$ is done by
the 3D cosine FFT algorithm; it takes $\mathcal{O}(N^{3}\log N)$ flops.
Operator $\mathcal{R}$ is defined by equations (\ref{E:diagonal3D}); it
consists in expanding the data for each face and for each time sample in the
2D cosine Fourier series (i.e. finding $g_{1,l,n}(t)$, $g_{2,k,n}(t),$ and
$g_{3,k,l}(t))$, in multiplying these functions by the scaled cut-off function
$\eta_{T},$ computing the 1D Fourier transform of the products at frequencies
$\omega_{k,l,n},$ and multiplying the results by some constant factors.
Expanding the data in the Fourier series takes $\mathcal{O}(N)\times
\mathcal{O}(N^{2}\log N)=\mathcal{O}(N^{3}\log N)$ flops. The only remaining
non-trivial operation here is computing the 1D Fourier transforms
$\widehat{\eta_{T}g_{...}}$ of products $\eta_{T}g_{...}$ at non-equally
spaced frequencies. In our simulation we simply applied the regular FFT to
compute $\widehat{\eta_{T}g_{...}}$ on a uniform grid, and interpolated using
third degree Lagrange polynomials, to obtain the values at the frequencies
$\omega_{k,l,n}$. The computational expense of this step is $\mathcal{O}%
(N^{3}\log N)$ flops. A more sophisticated approach to finding these values is
to use the proper version of the NSFFT \cite{Rokhlin}. This would yield more
accurate results at the price of longer computing time, although the
complexity would still remain $\mathcal{O}(N^{3}\log N)$ flops for this step.
We found, however, that the simple polynomial interpolation described above is
sufficiently accurate for the purposes of the present work.

To summarize, one iteration of the proposed method requires $\mathcal{O}%
(N^{3}\log N)$ flops, i.e. the algorithm is indeed asymptotically fast. The
total number of iterations needed to attain the convergence was from 1 to 4 in
our simulations, so, it is fair to say that the whole technique is implemented
in $\mathcal{O}(N^{3}\log N)$ flops.

\noindent\begin{figure}[t]
\begin{center}
\subfigure[]{\includegraphics[width=1.9in,height=1.9in]{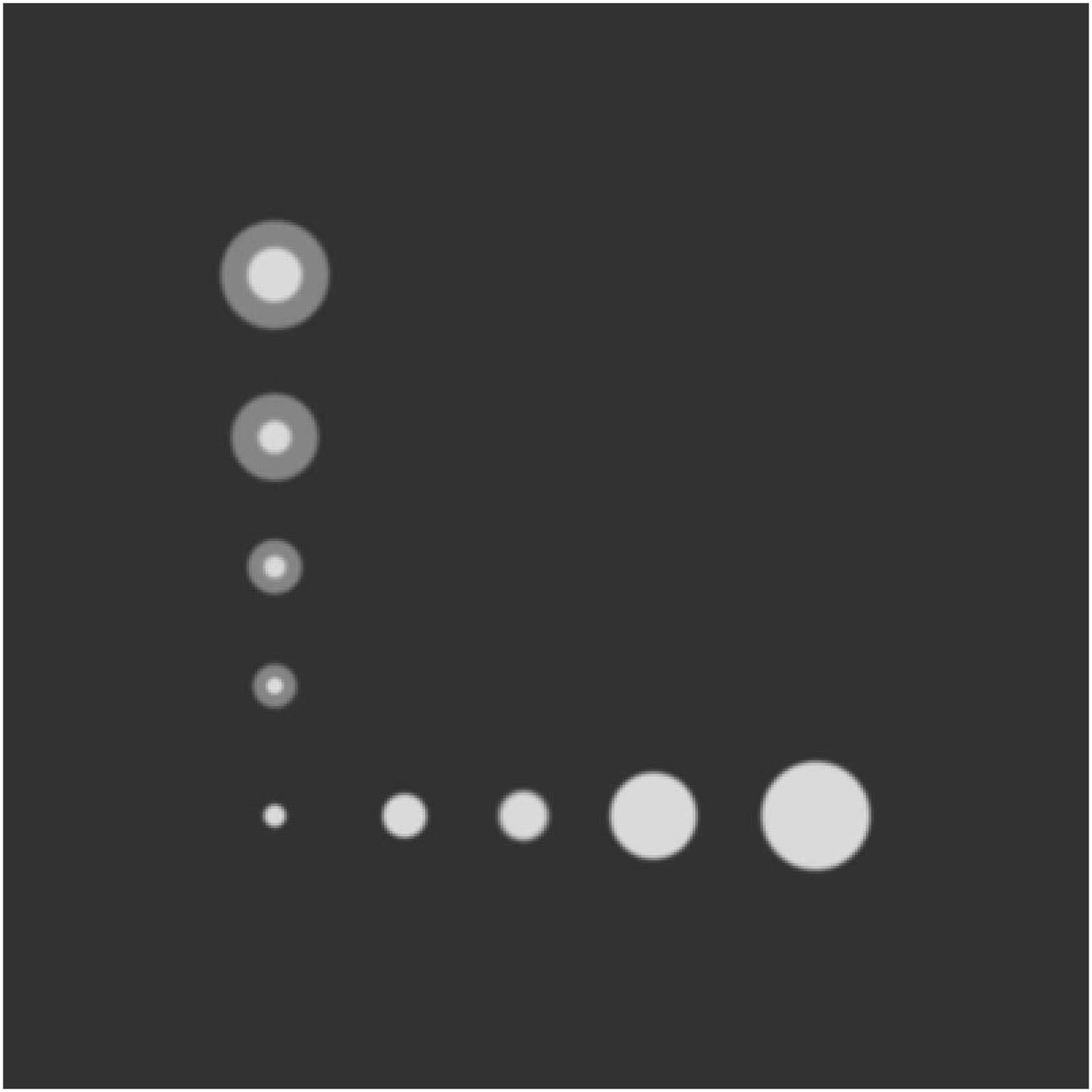}}
\subfigure[]{\includegraphics[width=1.9in,height=1.9in]{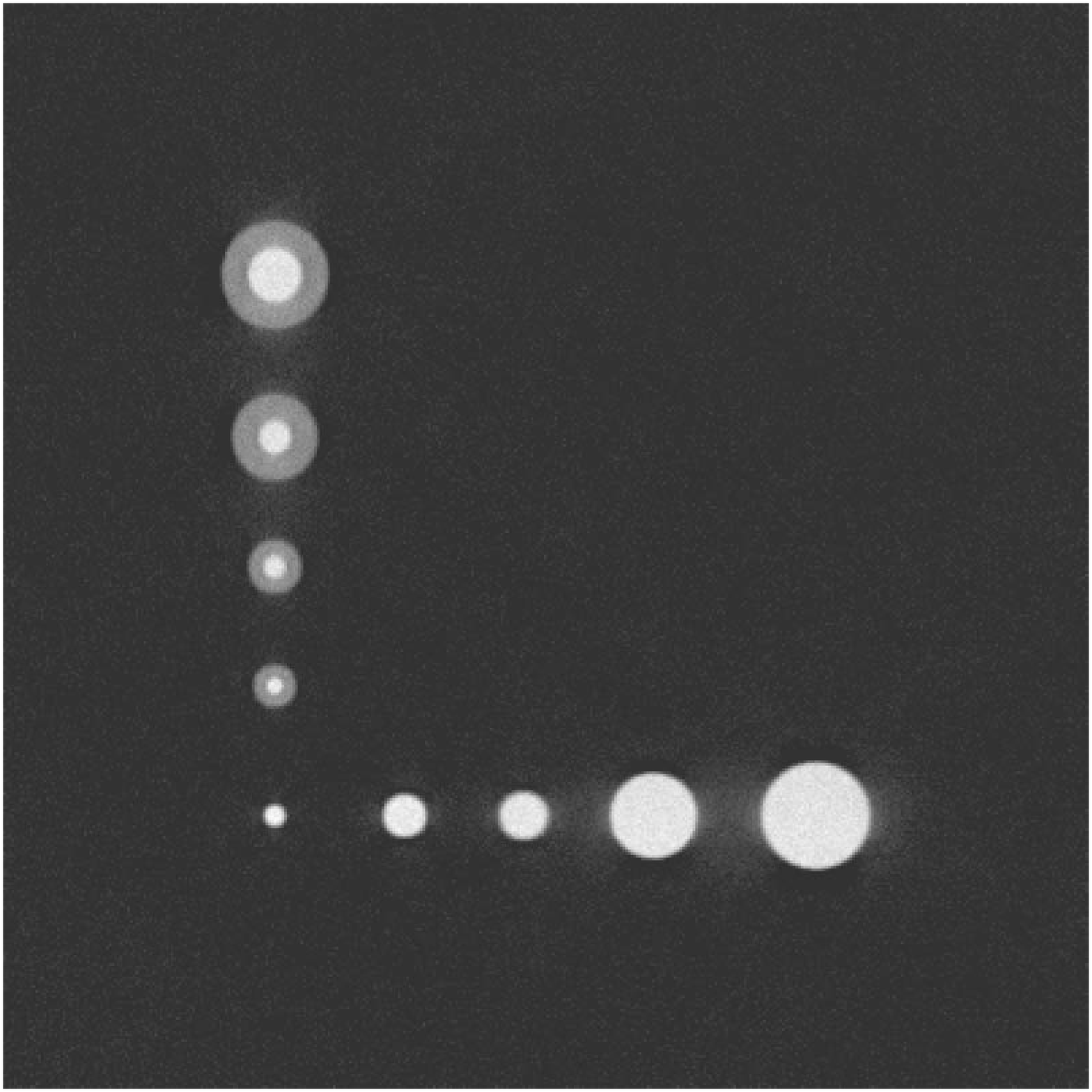}}
\subfigure[]{\includegraphics[width=1.9in,height=1.9in]{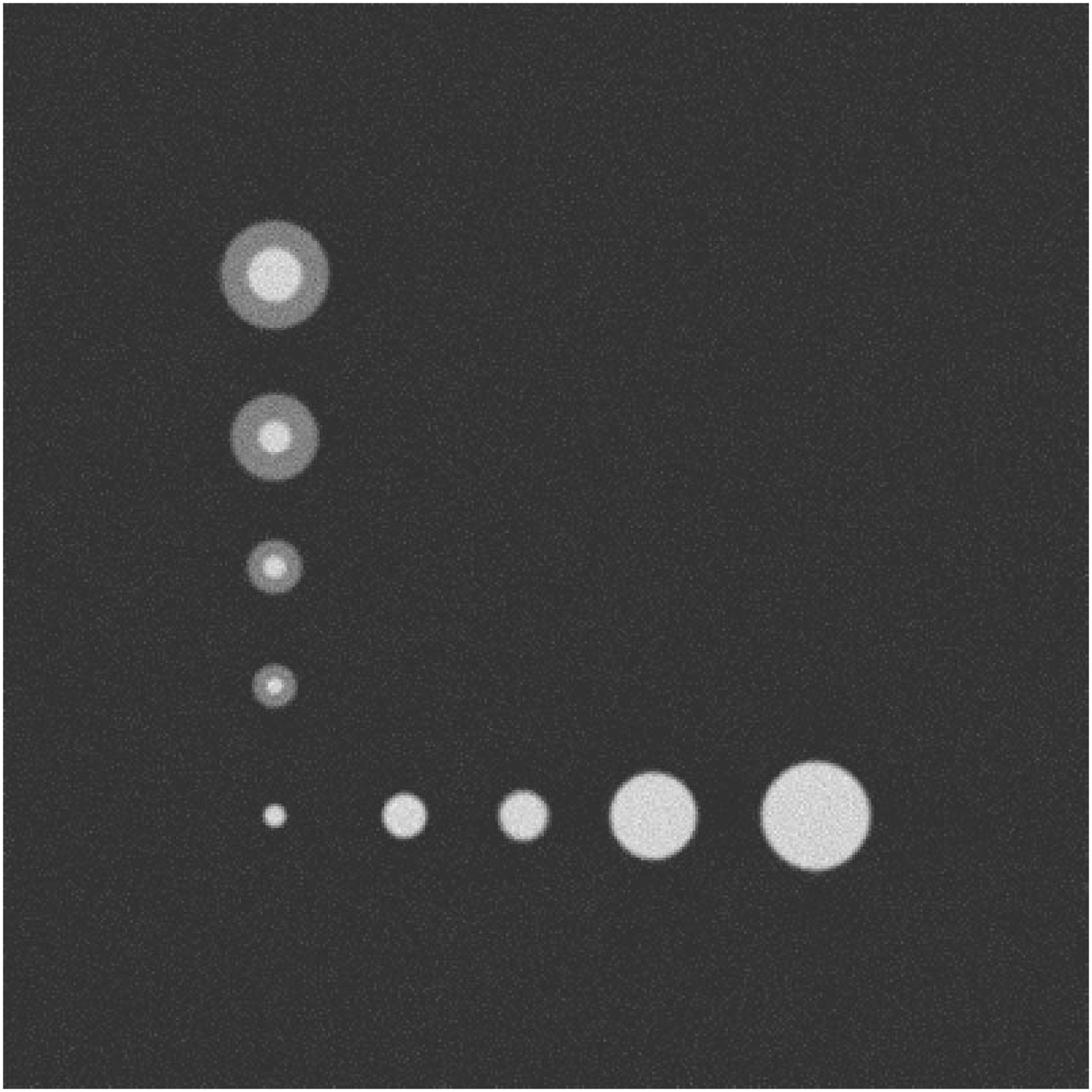}}
\end{center}
\caption{Simulation with $100\%$ noise in the data (a) Phantom (b) Initial
approximation $f^{(0)}(x)$ (c) 1-st iteration $f^{(1)}(x)$ }%
\label{F:n4}%
\end{figure}

\subsection{Simulations}

Performance of the proposed reconstruction algorithm was tested in numerical
simulations. The preliminary results in 2D were reported in \cite{CHK}. Below
we present 3D simulations.

\subsubsection{Measurements done on three faces}

As a numerical phantom representing the initial pressure $f(x)$ we used a
linear combination of several slightly smoothed characteristic functions of
balls of various radii, similar to the phantom utilized in \cite{Kun-MAET}.
The centers of the balls were chosen to lie on pair-wise intersections of the planes
$x_{j}=0.25,$ $j=1,2,3$. The cross section of the phantom by the plane
$x_{3}=0.25$ is shown in Figures \ref{F:n1}(a) and \ref{F:n4}(a). (If needed,
additional cross-sections can be found in \cite{Kun-MAET}.). The measurement
time $T$ in the two simulations described below was equal to $2$,  which
corresponds to the sound wave covering the shortest distance between the
parallel faces twice. For comparison, the measurement time in the standard
problem of TAT/PAT in 3D free space for a cube of this size would equal
$\sqrt{3}$.

The size of the discretization grid was $401\times401\times401$, i.e.
the reconstructed images contained about 64 million of unknowns. For simplicity
the time step was chosen to equal the spatial step (which was reasonable since
the model speed of sound was equal to 1). Therefore, $800$ time steps of the
forward problem were simulated. In both simulations the data were `measured'
on three mutually adjacent faces of the cubical domain.

Slices by the plane $x_{3}=0.25$ of the 3D images reconstructed from accurate
data are shown in Figure \ref{F:n1}, parts (b) and (c). Part (b) corresponds
to the initial approximation $f^{(0)}$; parts (c) represents the second
iteration $f^{(2)}$. The gray scale image in part (b) looks quite close to the
slice of the phantom shown in part (a). However, let us look at the Figure
(\ref{F:n2}) which shows the profiles of the functions along line $x_{2}
=x_{3}=0.25$. On can see that the first approximation $f^{(0)}$ (shown by the
dashed line) is not quite accurate. On the other hand, the profile
corresponding to the second iteration $f^{(2)}$ (represented by the gray solid
line) practically coincides with that of $f$ (black line).

In order to test the sensitivity of our reconstruction technique to the noise
always contained in real data, we added a strong noise component to the
simulated data. The noise was modelled by a normally distributed random
variable; it was scaled so that the noise intensity (as measured by the
standard $L_{2}$ norm) coincided with the intensity of the unperturbed signal.
Figure (\ref{F:n3}) demonstrates visually the high level of noise in the data;
the black line shows exact signal $g_{1}(t,0.5,0.5)$, while the gray line
represents the same signal with added noise.

It should be noted that most tomography modalities tend to amplify noise, and
in the presence of $100\%$ noise either the reconstructed image would be
overwhelmed by noise-related artifacts, or significant blurring would occur
due to the low-frequency filtration needed to regularize the reconstruction.
This does not happen here, since the singularities of the solution $p(t,x)$ do not
get smoothed on their way to the detectors due to the properties of the wave
equation. (Such low sensitivity to noise was observed previously
\cite{KuKuAET1,Kun-MAET} in other inverse problems related to PAT/TAT.)

Figure \ref{F:n4}, parts (b) and (c), contains images reconstructed from the
noisy data (the initial approximation, $f^{(0)}$, and the first iteration,
$f^{(1)}$, respectively). The noise is almost unnoticeable in the grayscale
images.
In addition to the low noise sensitivity of the method this phenomenon
is partially explained by the nature of the phantom (the volume of the
support of $f$ is actually quite small compared to the volume of the cube).

\begin{figure}[t]
\begin{center}
\includegraphics[width=6in,height=1.3in]{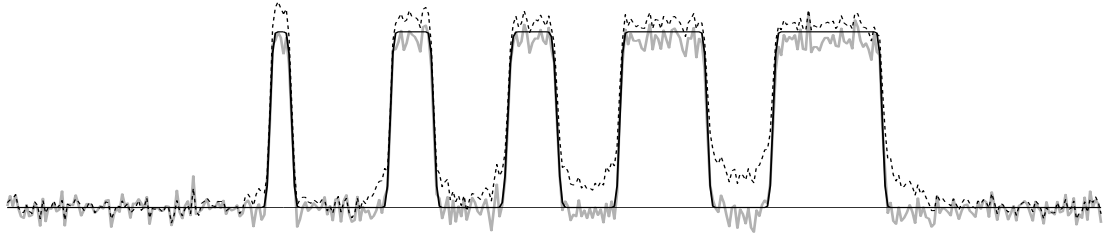}
\end{center}
\caption{Cross section of images in Figure~\ref{F:n4} along the line
$x_{2}=0.25$, $x_{3}=0.25$. The solid line represents the phantom, the dashed
line shows $f^{(0)}(x)$ and the gray line represents $f^{(1)}(x)$.}%
\label{F:n5}%
\end{figure}

Figure \ref{F:n5} presents the profiles of the cross sections by the line
$x_{2}=x_{3}=0.25$. One can see that a noticeable noise is present in the
reconstructions, and that the first iteration (gray line) does represent a
considerable improvement comparing to the initial approximation $f^{(0)}$. In
our simulation, in the presence of the strong noise the consecutive
approximations $f^{(j)},$ $j=2,3,4,..$. did not show any improvement over
$f^{(1)};$ moreover, a very slow growth in the level of noise was observed.

The computation time in the above simulations was about 9 minutes per
iteration (excluding the input/output time) on a 2.6 GHz processor of a
desktop computer. The code was not highly optimized, and it ran as
a single-thread (no parallelization).
Most of the elapsed time (7 min.) was spent computing the
forward problem (i.e. $\mathcal{W}f^{(K-1)}$). This is due to a large constant
factor implicit in the operation count of the NUFFT used on this step.

The computing time required by our algorithm is of the same order of magnitude as
that of the fast algorithm for the free-space space problem with a cube
surface acquisition \cite{Kun-ser} (estimated 4 to 5 min. on the
$401\times401\times401$ grid). Somewhat faster reconstruction time (about 1
min. for a slightly larger grid) was reported in \cite{Kun-fast} for a fast
algorithm for the free-space problem involving integrating linear detectors.
(The latter two problems are somewhat simpler computationally; in particular,
the algorithms are non-iterative). It should be noted that the
 time reversal algorithm for the free space problem implemented using finite differences
 (with complexity $\mathcal{O}(N^{4})$ flops) would take an estimated 3.5 hours on a grid of our
size, and methods based on a straightforward discretization of explicit
inversion formulas (whose complexity equals $\mathcal{O}(N^{5})$ flops) would
take several days to complete one reconstruction.

\subsubsection{Measurements done from all six faces}

We have also conducted simulations aimed at estimating the impact of
additional measurements (done on additional faces of the cube). Clearly, if
the data were measured on the other three faces of the cube (corresponding to
planes $x_{j}=1,$ $j=1,2,3)$, an algorithm very similar to described above
could be used for reconstruction. When the data are measured on all six faces,
the image is obtained simply by averaging the results obtained from the two
sets of three-face measurements. Our simulations demonstrate that with the
additional data the measurement time $T$ can be decreased. With measurement
time $T=1$ and a six-face acquisition the results of reconstruction were quite
similar to those corresponding to the three-face acquisition and $T=2$, as
described at the beginning of this section.

\section*{Acknowledgments}

The work of B.T. Cox was partially supported by the Engineering and Physical Sciences
Research Council, UK, and the work of L. Kunyansky was partially supported by
the NSF\ grants DMS-1211521 and DMS-0908243.



\end{document}